\documentclass{ifacconf}
\usepackage{amsmath, amssymb}
\usepackage{mathtools}
\usepackage{xcolor}
 \usepackage{natbib}
\usepackage{algorithm}
\usepackage{algpseudocode}
\newtheorem{theorem}{Theorem}

\newtheorem{Assumption}{Assumption}

\newtheorem{Definition}{Definition}
\newtheorem{Remark}{Remark}
\newtheorem{Example}{Example}
\usepackage[T1]{fontenc}

\begin{document}
\begin{frontmatter}

\title{Sum-of-Squares Certificates for\\ Almost-Sure Reachability of Stochastic Polynomial Systems} 

\author[first]{Arash Bahari Kordabad} 
\author[first]{Rupak Majumdar} 
\author[first]{Sadegh Soudjani}

\address[first]{Max Planck Institute for Software Systems, Kaiserslautern, Germany (e-mails: {\{arashbk,rupak,sadegh\}@mpi-sws.org})}

\begin{abstract}
In this paper, we present a computational approach to certify almost sure reachability for discrete-time polynomial stochastic systems by turning drift–variant criteria into sum-of-squares (SOS) programs solved with standard semidefinite solvers. Specifically, we provide an SOS method based on two complementary certificates: (i) a drift certificate that enforces a radially unbounded function to be non-increasing in expectation outside a compact set of states; and (ii) a variant certificate that guarantees a one-step decrease with positive probability and ensures the target contains its nonpositive sublevel set. We transform these conditions to SOS constraints. For the variant condition, we enforce a robust decrease over a parameterized disturbance ball with nonzero probability and encode the constraints via an S-procedure with polynomial multipliers. The resulting bilinearities are handled by an alternating scheme that alternates between optimizing multipliers and updating the variant and radius until a positive slack is obtained. Two  case studies illustrate the workflow and certifies almost-sure reachability.
\end{abstract}


\end{frontmatter}

\section{Introduction}
\label{sec:intro}

Ensuring that the trajectories of a dynamical system eventually reach a given target set, known as reachability, is a fundamental objective in control and formal verification~\citep{baier2008principles}. In deterministic settings, classical constructs such as Lyapunov and barrier functions provide crisp, certificate-based guarantees without the need for having closed-form expressions of the trajectories~\citep{ames2019control}. However, once uncertainty enters the dynamics, these certificates must be adapted to account for statistics of the uncertainty~\citep{kordabad2024control,lavaei2022automated}. 

For stochastic systems, the computation of reachability probability hinges on solving the almost sure reachability problem, which is characterizing systems that a target set can be reached with probability one \citep{junges2021enforcing}. Recent work has established necessary and sufficient conditions for almost-sure reachability of general discrete-time stochastic systems via drift and variant functions with one-step ahead conditions~\citep{majumdar2024necessary}.  The drift is radially unbounded and does not increase in expectation outside a compact set. This ensures that trajectories of the system do not escape to infinity almost surely~\citep{meyn2012markov}. On the other hand, the variant is positive outside the target set and exhibits a uniform one-step decrease with positive probability on suitable sublevel sets of the drift. This provides a positive probability for the trajectories to move towards the target set. The drift and variant conditions together are shown to be necessary and sufficient for almost sure reachability. For linear systems with additive disturbance, a full characterization of almost sure reachability has been provided in~\citep{kordabad2025certificates} based on the system matrices using these drift–variant conditions. However, their algorithmic computation for nonlinear systems remains challenging due to the conditions on probability and expectation of the certificates with respect to the randomness over state space.

Sum-of-squares (SOS) programming certifies nonnegativity of polynomials by representing them as sums of squares, thereby converting the search for a certificate into a convex semidefinite program (SDP) solvable by off-the-shelf numerical solvers. SOS has become a standard tool in control and verification e.g., for synthesizing Lyapunov and barrier certificates and enforcing safety and performance constraints for polynomial systems~\citep{papachristodoulou2005tutorial,schneeberger2023sos,jagtap2020formal}. Consequently, conditions such as radial unboundedness, nonpositive drift outside a compact set, target containment, and robust one-step decrease can be cast as SOS programs; and, as in this work, drift–variant conditions will be evaluated using SOS to certify almost-sure reachability.

\cite{kordabad2025certificates} have shown that restricting the search for certificates to fixed templates can compromise completeness, and the focus there is on characterizing almost-sure reachability for linear systems. Using variant functions and supermartingales,~\cite{majumdar2025sound} provide sound and complete proof rules for qualitative and quantitative termination of probabilistic imperative programs with demonic bounded nondeterminism. From a computational standpoint, however, constructing such functions remains challenging: unconstrained searches are typically intractable, and while restricting templates may lose completeness guarantees~\citep{ahmadi2011globally}, it enables a practical, solver-ready framework. This motivates our choice to work within a specific, analyzable class and to develop implementable tools that, despite conservatism, make the certificates computational and usable in practice.

In this paper, we focus on the computational side of the almost sure reachability certificates and show how to turn the drift–variant criteria into tractable optimization for polynomial stochastic systems. Our approach uses SOS programming to encode the drift and variant requirements as convex semidefinite constraints. Concretely, we select polynomial templates for the drift and variant and: (i) enforce radial unboundedness and non-increasing drift in expectation outside a compact set via SOS; (ii) reformulate the probabilistic variant condition by introducing a parameterized disturbance subset, as a ball, and imposing a robust one-step decrease for all disturbances in that ball, while treating the ball’s radius as a decision variable to ensure the subset has nonzero probability; and (iii) encode target containment and other constraints by SOS multipliers. The resulting SOS program contains bilinear couplings between the variant, multipliers and other scalars. Using a similar approach to, e.g.,~\cite{yin2021backward,lin2025modified}, we resolve these with an alternating scheme that fixes the variant and radius to optimize multipliers and margins, then fixes multipliers to update the variant, until a strictly positive slack is obtained. This yields a practical, solver-ready pipeline that leverages standard SDP solvers and the well-developed SOS toolbox. 

\textbf{Contribution:} The main contributions are threefold: (i) a precise SOS encoding of the drift–variant  criteria for polynomial systems that separates structural constraints from design choices; (ii) a robust realization of the probabilistic variant via an optimizable disturbance ball with guaranteed positive probability, together with target-containment multipliers; and (iii) an alternating synthesis routine that resolves the induced bilinearities and returns certificates with quantitative margins and a certified radius. Two nonlinear polynomial examples demonstrate the approach.

\textbf{Outline:} Section~\ref{sec:preliminaries} recalls preliminaries on polynomial stochastic systems and the drift–variant theorem for almost-sure reachability. Section~\ref{sec:SOS} casts the drift and variant conditions as SOS constraints and presents the corresponding algorithms, including the alternating synthesis scheme. Section~\ref{sec:illustrative} demonstrates the approach on two nonlinear stochastic systems. Finally, Section~\ref{sec:conclusion} gives discussions and concluding remarks on the approach of the paper.

\textbf{Notation:} We write $\mathbb{R}$ for the reals, $\mathbb{R}_{\ge 0}$ for the nonnegative reals, $\mathbb{R}_{> 0}$ for the positive reals, and $\mathbb{N}_{\ge 0}$ for the nonnegative integers. Vectors are columns and $\|x\|$ is the Euclidean norm. Expectation and probability are denoted by $\mathbb{E}[\cdot]$ and $\mathbb{P}(\cdot)$. For multivariate polynomials, $\mathbb{R}[x]$ is the ring in variables $x$ with real coefficients, $\deg (\cdot)$ is the total degree, and $\Sigma[x]$ is the SOS cone; analogously, $\Sigma[x,w]$ is the SOS cone in joint variables $(x,w)$. A set \( S \subseteq \mathbb{R}^n \) is called {open} if for every point \( x \in S \), there exists a neighborhood of \( x \) entirely contained in \( S \). A set is {closed} if its complement is open. A set $S \subset \mathbb{R}^n$ is {bounded} if there exists an $M\in \mathbb{R}_{\geq 0}$ such that $\|x\| \leq M$ for all $x \in S$. A set is said to be {compact} if it is closed and bounded. 

\section{Preliminaries and Problem Statement}
\label{sec:preliminaries}
In this section, we first provide the preliminaries on sum-of-squares (SOS) polynomials. We then introduce the discrete-time polynomial stochastic system (dt-PSS) model considered in this paper, followed by the formulation of the almost-sure reachability problem. Finally, we recall the recently developed drift–variant conditions, which provides necessary and sufficient conditions for almost-sure reachability and forms the basis of our SOS formulation.
\begin{Definition}
For $x \in \mathbb{R}^n$, a multivariate polynomial $p(x)\in \mathbb{R}[x]$ is an SOS, and denoted as $p(x)\in \Sigma[x]$, if there exist some polynomials $f_i(x)\in\mathbb{R}[x]$, $i = 1, \ldots, M$, such that
$
    p(x) = \sum_{i=1}^{M} f_i^2(x)
$.
\end{Definition}
An equivalent characterization of SOS polynomials has been provided by~\cite{parrilo2000structured}. Specifically, a  polynomial $p(x)$ of degree $2d$ is SOS if and only if there exists $Q \succeq 0$ and a vector of monomials $z(x)$ containing all monomials in $x$ of degree $\leq d$ such that
$
    p(x) = z(x)^\top Q z(x). $
Note that all SOS polynomials are nonnegative, but the converse is not necessarily true; that is, there exist nonnegative polynomials that are not SOS (see e.g.,~\cite{ahmadi2012convex}). However, this characterization enables SOS problems to be formulated as semidefinite programs (SDPs), which can be efficiently solved using convex optimization techniques.

\smallskip
We consider an underlying dynamical model given by a discrete-time polynomial stochastic system (dt-PSS), defined as a tuple $\mathfrak{S}=(X,W,w, f)$, where $X\subseteq\mathbb{R}^n$ is the state space of the system,
$W\subseteq\mathbb{R}^m$ is the disturbance space. The disturbance is a sequence
$w:=\{w_k:\Omega\rightarrow W,k\in\mathbb{N}_{\geq 0}\}$ of independent and identically distributed (i.i.d.) random variables defined on a sample space $\Omega$. Each $w_k$ takes values in $W$ and is sampled according to a probability measure $\mathbb{P}_w$. The polynomial map $f:X\times W\rightarrow X$  characterizes the state evolution of the system according to
	\begin{equation}
	x_{k+1}=f(x_k,w_k),
	\label{eq:state_evolution}
	\end{equation}
where $x_k\in X$ is the state at time $k\in\mathbb{N}_{\geq 0}$, $x_0$ is an initial condition and $w_k\in W$ is a random variable that satisfies the following assumption. 

\begin{Assumption}
\label{ass:disturbance}
The disturbance sequence $\{w_k\}_{k\in\mathbb{N}_{\ge 0}}$ consists of random variables $w_k \in W \subseteq \mathbb{R}^m$ with probability measure $\mathbb{P}_w$. 
Each $w_k$ is assumed to have zero mean, $\mathbb{E}[w_k]=0$, without loss of generality, and its (possibly unbounded) support contains an open ball centered at the origin. Moreover, the distribution $\mathbb{P}_w$ admits finite moments up to the degree required by the polynomial expressions appearing in the approach. 
\end{Assumption}

In this paper, the function $f$ is assumed to be a polynomial in both arguments. More specifically, we suppose that $f = (f_1, \dots, f_n)$ is a vector of polynomials and each coordinate function $f_i(x, w)\in \mathbb{R}[x,w]$ is a multivariate polynomial.

We now state the almost sure reachability question addressed in this work.

\textbf{Almost Sure Reachability Problem.}
\label{prob:asr}
Given a dt-PSS $\mathfrak{S}$ and a bounded open target set $G \subseteq X$, determine whether for all initial conditions $x_0 \in X$, the system trajectory $\{x_k\}_{k=0}^\infty$ hits the target set $G$ almost surely, i.e., 
\begin{equation}\label{eq:reach}
    \forall x_0\in X,\quad \mathbb{P}\left( \exists\, k \in \mathbb{N}_{\geq 0} : x_k \in G  \right) = 1.
\end{equation} 

To characterize almost sure reachability, we adopt the framework proposed by~\cite{majumdar2024necessary}, which provides a pair of conditions that are necessary and sufficient for ensuring that the trajectories reach a given target set \( G \) with probability one.

The theorem by~\cite{majumdar2024necessary} applies to more general discrete-time stochastic systems that are \emph{weak Feller}. A stochastic system is said to be weak Feller if its transition maps bounded continuous functions to continuous functions, ensuring continuity of the expected next-state value with respect to the current state. For a precise definition of the weak Feller property see e.g.,~\cite{meyn2012markov}. Note that one can readily verify that the dt-PSS $\mathfrak{S}$ defined above is weak Feller. We recall these conditions next.

\smallskip

\noindent\textbf{V1: Drift Criterion.} There exists a drift function $V:X\rightarrow\mathbb{R}_{\geq 0}$ with $\lim_{\|x\|\rightarrow\infty}V(x)=\infty$ and a compact set $C\subseteq X$ satisfying
\begin{equation}
\label{eq:drift}
    \Delta V(x) :=\mathbb{E}\left[ V(f(x, w))\,|\,x \right] - V(x) \le 0,\quad \forall x\in C^c,
\end{equation}
where $C^c$ is the complement of $C$.

Intuitively, $V$ acts as a Lyapunov or supermartingale function: outside a (possibly large) compact set, $V$ does not increase in conditional expectation, which prevents trajectories from escaping to infinity with positive probability, which is necessary for satisfying~\eqref{eq:reach}. 

\smallskip

\noindent\textbf{V2: Variant Criterion.} For a function $V$ satisfying the drift criterion \textbf{V1}, there exists a function $U : X \to \mathbb{R}$ called the \emph{variant}, and three supporting functions $H: \mathbb{R}_{>0} \to \mathbb{R}$, $\delta : \mathbb{R}_{>0} \to \mathbb{R}_{>0}$, and $\epsilon : \mathbb{R}_{>0} \to \mathbb{R}_{>0}$ such that for all $r\in \mathbb{R}_{>0}$ and $x \in X$, the implication 
$
    V(x) \leq r \implies U(x) \leq H(r)
$
holds, and 
\begin{equation}\label{eq:variant}
    \mathbb{P}_w(U(f(x,w))-U(x)\le -\delta(r))\ge \epsilon(r),
\end{equation}
for all $x$ satisfying $V(x) \leq r$ and $U(x)> 0$.

Thus, on each drift sublevel set $\{x: V(x)\le r\}$, the variant $U$ is bounded above by $H(r)$ and makes one-step progress with a strictly positive probability $\epsilon(r)$ by at least a strictly positive decrement $\delta(r)$ whenever $U(x)>0$. Next, we state the almost-sure reachability theorem using \textbf{V1-V2}.

\begin{theorem}[Almost sure reachability]
\label{thm:nece_suff}
    For a dt-PSS $\mathfrak{S}$ and open bounded target set $G$, if there exists a function $V$ satisfying criterion \textbf{V1} and a variant $U$ satisfying criterion \textbf{V2} associated with $V$, such that
    \begin{equation}\label{eq:G:U}
     G \supset \{x \in X \mid U(x) \leq 0\},  
    \end{equation}
then~\eqref{eq:reach} holds. Moreover, if \eqref{eq:reach} holds, then there are $V$ and $U$ that satisfy \textbf{V1} and \textbf{V2}, and such that~\eqref{eq:G:U} holds.
\end{theorem}
\textit{Proof.} See~\cite{majumdar2024necessary}. \hfill $\blacksquare$

We provide the following polynomial system to demonstrate two points:
(a) The variant condition in \textbf{V2} is not sufficient by itself to guarantee almost sure reachability. The drift condition \textbf{V1} is also needed to ensure the solution does not diverge to infinity.
(b) Assuming constant $\delta$ and $\epsilon$ in \textbf{V2} reduces the chance of satisfying this condition if the state space of the system is unbounded.
\begin{Example}\label{example1}
    Let us consider the following one-dimensional polynomial system and the target set \begin{equation*}
        x_{k+1} = x_k + x_k^2 w_k, \quad w_k \sim \mathcal{U}[-1,1], \quad G=(-2,2).
    \end{equation*} 
 We first show that this system diverges to infinity with positive probability, thus it does not satisfy the condition \eqref{eq:drift} and is not almost sure reachable. Assume the initial condition satisfies $|x_0| > 4$. At each time step $k$, we select a subset of the disturbance space defined as
 $$W_k=\{w_k\in[-1,1] \,\mid\,|1+x_kw_k|>2\}.$$
 This implies $|x_{k+1}|=|x_k||1+x_kw_k|>2|x_k|$ for all $w_k\in W_k$, therefore $|x_k|\geq 2^{k+2}$ for all $k\geq 0$. Moreover, this choice of disturbance sets yields $w_k>\frac{1}{x_k}$ or $w_k<-\frac{3}{x_k}$ for $x_k>0$ and $w_k<\frac{1}{x_k}$ or $w_k>-\frac{3}{x_k}$ for $x_k<0$. In both cases,  the probability of $w_k\in W_k$ is $1-\frac{2}{|x_k|}$, which is lower bounded by $1-\frac{1}{2^{k+1}}$. Hence, the probability of obtaining such a diverging trajectory is bounded below by $p=\Pi_{k=0}^\infty( 1-\frac{1}{2^{k+1}})$. One can observe that, $\ln(p)=\sum_{k=0}^\infty\ln( 1-\frac{1}{2^{k+1}})$. Using the inequality $\ln(\eta)\geq -\frac{\eta}{1-\eta}\geq -2\eta$ for all $\eta\in(0, \frac{1}{2}]$, we obtain $\ln(p)\geq -\sum_{k=0}^\infty 2^{-k}=-2$. Therefore, $p\geq \mathrm{e}^{-2}>0$. We then use results from~\cite{meyn2012markov} to conclude that no drift function can satisfy condition~\textbf{V1}. Recall that a system is said to be \emph{non-evanescent} if, for all initial conditions $x_0$, we have $\mathbb{P}\left(\|x_k\|\to \infty \right)=0$. However, we showed that the system violates the non-evanescence condition. Therefore, by Theorem 9.4.1 by~\cite{meyn2012markov}, it follows that there exists no function $V$ and a compact set $C$ satisfying the drift condition~\eqref{eq:drift}.

We now show that the candidate variant function $U(x)=x^2-1$ satisfies condition~\textbf{V2}. First, note that $\{x\mid U(x)\leq 0\}\subset G$. Let us assume that a constant $\delta(r)\equiv\delta<1$ satisfies~\eqref{eq:variant}. Then, the inequality inside the probability in~\eqref{eq:variant} can be rewritten as
\begin{align}\label{eq:exam1:ineq}
    &(x+x^2w)^2-x^2\leq -\delta \Leftrightarrow (x+x^2w)^2\leq x^2-\delta \nonumber\\
    &|x+x^2w|\!\leq \!\sqrt{x^2-\delta}\! \Leftrightarrow\! -\sqrt{x^2-\delta}\leq  x+x^2w \leq \sqrt{x^2-\delta}\nonumber\\
    &\Leftrightarrow-\frac{1}{x}-\frac{\sqrt{x^2-\delta}}{x^2}\leq  w \leq -\frac{1}{x}+\frac{\sqrt{x^2-\delta}}{x^2},
\end{align} 
for all $x$ with $x^2> 1$ (i.e., $U(x)>0$). Note that $-1\leq -\frac{1}{x}-\frac{\sqrt{x^2-\delta}}{x^2}$ and $-\frac{1}{x}+\frac{\sqrt{x^2-\delta}}{x^2}\leq 1$ for all $x\in G^c$. The probability that inequality~\eqref{eq:exam1:ineq} holds is therefore
$
   \bar \epsilon(x)=\frac{\sqrt{x^2-\delta}}{x^2}>0,
$ 
and it satisfies $\lim_{|x|\rightarrow\infty}\bar \epsilon(x)=0$. Hence,
$
  \epsilon (r)=\inf_{V(x)\le r}   \bar \epsilon(x),
$
is a valid $\epsilon (r)$ for~\eqref{eq:variant} with $ \lim_{r\rightarrow \infty}\epsilon (r)=0$. Therefore, no constant $\epsilon$ and $\delta$ can satisfy condition~\textbf{V2} for this choice of variant function, which is a valid variant function.
\end{Example}

Verifying \textbf{V1}–\textbf{V2} for given $(V,U)$ only requires the one-step map $f$ and the distribution of $w$, without reasoning over infinite trajectories. The theorem is, however, \emph{existential}, i.e., it asserts that suitable $(V,U)$ are necessary and sufficient for almost-sure reachability, but it does not prescribe how to construct them or in which function classes to search.

This paper focuses on the {computational aspects} of the above certificates, since the practical value of \textbf{V1}–\textbf{V2} hinges on solving them efficiently. For polynomial systems, the SOS framework and its reduction to SDPs enable computation with off-the-shelf solvers. By selecting polynomial templates for $V$ and $U$, we obtain SOS programs that certify the required conditions on the state space.  Note that Theorem~\ref{thm:nece_suff} asserts that the drift-variant criteria are, in principle, necessary and sufficient for almost-sure reachability. However, once we fix templates for the certificates such as polynomial functions, the resulting SOS programs provide only a sufficient computational test: if a valid certificates exist outside the chosen class, the SOS search may fail to find it. See, e.g., the discussion and examples in~\cite{kordabad2025certificates}, which illustrate that fixing polynomial templates can imply incompleteness even for polynomial systems. This loss of completeness is the trade-off for convex tractability and scalable tooling.


\section{Sum-of-Squares Certificates}\label{sec:SOS}
In this section, we turn the theoretical criteria \textbf{V1} (drift) and \textbf{V2} (variant) into concrete SOS programs. We reformulate the almost-sure reachability conditions as SOS constraints once polynomial templates for $V$ and $U$ are fixed. The resulting constraints (i) enforce the drift inequality~\eqref{eq:drift} in expectation outside a compact set, (ii) certify a one-step decrease of $U$ on a disturbance set of positive probability as expressed in~\eqref{eq:variant}, and (iii) guarantee target containment~\eqref{eq:G:U}. 

\subsection{SOS-based Drift function}
We search for a polynomial drift function $V(x)$ in \textbf{V1} within an SOS template, i.e., $V(x)\in \Sigma[x]$. Then it satisfies $V(x)\geq 0$ by construction.
To enforce  the radially unboundedness of $V$, we impose 
\begin{equation}\label{eq:unbound:V}
     V(x) -\gamma_0 x^\top x + \lambda_0 \in \Sigma[x],
\end{equation}
with decision scalars $\gamma_0>0$ and  $\lambda_0\in\mathbb{R}$. This ensures that $V(x) \geq  \gamma_0 x^\top x - \lambda_0$ and hence $V(x)\to \infty$ as $\|x\| \to \infty$.
    
To satisfy \eqref{eq:drift}, first observe that since both \(f\) and \(V\) are polynomials, and assuming the disturbance admits finite moments up to a sufficiently large order, the composition \(V(f(x,w))\) is itself a polynomial in \((x,w)\).

\begin{Remark}
The SOS-based formulation relies on the ability to compute or bound the moments of the disturbance distribution, since expectations of polynomial functions of the noise appear explicitly in the drift and variant constraints. The finite-moment assumption is satisfied in most practical settings, including zero-mean i.i.d.\ disturbances with known moments arising from Gaussian, uniform, or other light-tailed distributions supported on a bounded set such as a box or a ball. 
\end{Remark}

The drift condition~\eqref{eq:drift} can be certified by the following SOS constraint:
\begin{equation}\label{eq:deltaV}
    -\Delta V(x)- \gamma_1x^\top x+\lambda_1 \in \Sigma[x],
\end{equation}
for some $\gamma_1>0$ and $\lambda_1\in\mathbb{R}$. This condition guarantees that for all $x$ outside of the compact set $C=\{x\in\mathbb{R}^n\,|\,x^\top x\leq {\lambda_1}/{\gamma_1}\}$, we have $\Delta V(x)\leq 0$ and therefore~\eqref{eq:drift} is fulfilled. Note that $\lambda_1<0$ corresponds to the empty $C$. With the SOS parameterization of \(V\) and the constraints \eqref{eq:unbound:V}–\eqref{eq:deltaV}, feasibility of the SDP produces a radially unbounded polynomial \(V\) whose expected drift is nonpositive outside a compact set \(C\), thereby certifying the drift part \textbf{V1}. Algorithm~\ref{alg:drift-sos} summarizes the steps.

\begin{algorithm}[t]
\caption{SOS-Based Drift Function}
\label{alg:drift-sos}
\begin{algorithmic}[1]
\State \textbf{Input:} Polynomial dynamics $f(x,w)$, desired drift function degree and disturbance moments.
\State \textbf{Output:} Polynomial drift function $ V(x)$.
\State \textbf{// Step 1: Construct the drift function}
\State Define drift function  $V(x)\in \Sigma[x]$.
\State Form expected value: $\mathbb{E}_w[V(f(x,w))]$ using symbolic expansion of $V(f(x, w))$ and disturbance moments.
\State Compute $\Delta V(x) := \mathbb{E}_w[V(f(x, w))] - V(x)$.
\State \textbf{// Step 2: Enforce the radially unboundedness}
\State Introduce scalar decision variables $\gamma_0>0$,  $\lambda_0$ and impose \eqref{eq:unbound:V}.
\State \textbf{// Step 3: Enforce the drift condition}
\State Introduce scalar decision variables $\gamma_1>0$, $\lambda_1$ and impose~\eqref{eq:deltaV}.
\State \textbf{// Step 4: SOS solver}
\State Solve SOS constraints to get $V(x)$, $\gamma_0$, $\lambda_0$, $\gamma_1$, and $\lambda_1$.
\State \textbf{Return:} Drift function $V(x)$ if feasible solution is found.
\end{algorithmic}
\end{algorithm}

\subsection{SOS-based Variant Function}
The variant criterion \textbf{V2} is distinctive in that it requires a uniform one-step decrease with strictly positive probability over the state space. This type of probabilistic, globally quantified decrease is not extensively studied in the literature and is computationally challenging: it couples the geometry of the state space with the distributional support of the disturbance. Our SOS construction below turns this requirement into tractable constraints. We select a polynomial template $U(x)\in\mathbb{R}[x]$ and define, 
$
    H(r):=\sup_{V(x)\leq r} U(x).
$
Because \(V\) is polynomial and (by \textbf{V1}) radially unbounded, every sublevel set \(\{x: V(x)\le r\}\) is compact; since \(U\) is polynomial, \(H(r)\) is finite and attained. This guarantees the implication \(V(x)\le r \Rightarrow U(x)\le H(r)\) used in \textbf{V2}.

To satisfy~\eqref{eq:variant} without integrating over the full noise law, we carve out a measurable subset $W_0\subseteq W$ of non-zero probability, i.e., $\mathbb P(w\in W_0)\geq \epsilon,$ for some $\epsilon>0$, and enforce a robust decrease on \(W_0\). Specifically, we enforce a uniform one-step decrease on $W_{0}$, as follows:
\begin{equation}\label{eq:dec:W0}
   U(f(x,w))-U(x)\le -\delta,\quad \forall w\in W_0,
\end{equation}
for some $\delta>0$ and all $x$ with $U(x)>0$. Then \eqref{eq:variant} holds with $\delta(r)\equiv\delta$ and $\epsilon(r)\equiv\epsilon$. We take $W_0$ as the Euclidean ball
$
    W_{0}=\{w\in \mathbb{R}^m\,:\, w^\top w\leq \rho\},
$
with radius parameter $\rho>0$. Under our standing assumption that the support of \(w\) contains an open ball around the origin, there exists a \(\rho\) such that \(\mathbb{P}(w\in W_0)>0\); shrinking \(\rho\) if needed also guarantees \(W_0\subseteq W\).  A sufficient condition, based on S-procedure, for the decrease condition in \textbf{(V2)} is then the existence of a polynomial $U\in\mathbb{R}[x]$ and scalars $\delta > 0$ and $\rho > 0$, together with SOS multipliers $\Lambda(x,w)$, and $M(x,w)$ such that
\begin{subequations}\label{eq:U:SOS}
   \begin{align}
U(x) - U(f(x,w))& - \delta - \Lambda(x,w)\bigl(\rho - w^\top w \bigr) \nonumber\\ - M(x,w) U(x)&\in \Sigma[x,w],\label{eq:variant:SOS}\\ \Lambda(x,w)&\in \Sigma[x,w],\\ M(x,w)&\in \Sigma[x,w].
\end{align} 
\end{subequations}
For all $w$ with $w^\top w\leq \rho$ (i.e., $w\in W_0$), and for all $x$ with $U(x)>0$, the last two terms in~\eqref{eq:variant:SOS} are non-positive since $\Lambda(x,w)$, and $M(x,w)$ are nonnegative and this results in $U(x) - U(f(x,w)) - \delta\geq 0$ on \(W_0\) whenever \(U(x)>0\).

For the target set containment condition in~\eqref{eq:G:U}, suppose that $G$ is given by sub-level of some polynomial functions, as follows:
\[
G = \bigl\{\,x\in\mathbb{R}^n \;\big|\; g_i(x)< 0,\ i=1,\dots,I \,\bigr\},
\]
with $g_i\in\mathbb{R}[x]$. Using S-procedure for~\eqref{eq:G:U}, we require: 
\begin{subequations}\label{eq:GU:SOS}
\begin{align}
-g_i(x)+S_i(x)U(x)-\alpha_i &\in \Sigma [x],\label{eq:GU:SOS1}\\S_i(x)&\in \Sigma [x],\,
\end{align}
\end{subequations}
for some $\alpha_i>0$ and all $i=1,\ldots, I$. For all $x$ with $U(x)\le0$, we have $S_i(x)\,U(x)\le0$, then \eqref{eq:GU:SOS1} yields $-g_i(x)\ge\alpha_i>0$. Thus~\eqref{eq:G:U} holds.  We then solve the following optimization by gathering all the constraints and letting the scalars $(\rho, \alpha_{1:I},\delta)$ to be part of decision variables:
\begin{align}\label{eq:opt:U}
  &  \max_{U,S_{1:I},\Lambda,M,\rho, \alpha_{1:I},\delta,\varepsilon} \quad \varepsilon,\\
   &\qquad\quad\,\, \text{s.t.}\qquad\quad\quad\, \eqref{eq:U:SOS}, \eqref{eq:GU:SOS},\nonumber\\
   & \qquad\qquad\qquad\qquad\,\,\,\, \alpha_i\geq \varepsilon, \,\, i\!=\!1,\ldots, I,\nonumber \\ &\qquad\qquad\qquad\qquad\,\,\,\,\delta\geq \varepsilon,\quad \rho>0.\nonumber
\end{align}
We relax the positivity requirements $\alpha_i>0$ and $\delta>0$ by introducing a slack variable $\varepsilon$. The variable $\varepsilon$ relaxes the positivity margins early on and is driven to a strictly positive value. The problem in~\eqref{eq:opt:U} is convex in \((\Lambda,M,S_{1:I})\) when \((U,\rho)\) are fixed, and convex in \((U,\rho)\) when \((\Lambda,M,S_{1:I})\) are fixed, but not jointly convex due to bilinear couplings \((\Lambda,\rho)\) and \((S_i,U),(M,U)\). To handle the bilinear couplings in~\eqref{eq:opt:U}, similar to e.g.,~\cite{safe2022, lin2023secondary}, we adopt an alternating SOS program in which each subproblem is a convex SDP.  Starting from an initial variant $U_0$ and radius $\rho_0$, we iterate:
\begin{enumerate}
\item \textit{Multiplier step:} Set $(U, \rho)=(U_k, \rho_k)$ and maximize $\varepsilon$ over $(S,\Lambda,M, \alpha_{1:I},\delta,\varepsilon)$ subject to the constraints~\eqref{eq:opt:U} to obtain $(S_{1:I,k},\Lambda_k,M_k, \alpha_{1:I,2k-1},\delta_{2k-1},\bar \varepsilon_{1})$.
\item \textit{Variant step:} Fix $(S,\Lambda,M, \rho)=(S_{1:I,k},\Lambda_k,M_k, \rho_k)$ and maximize $\varepsilon$ over $(U,\alpha_{1:I},\delta,\varepsilon)$ subject to the constraints~\eqref{eq:opt:U} to obtain $(U_{k+1},\alpha_{1:I,2k},\delta_{2k},\bar \varepsilon_{2})$.
\item \textit{Update and stopping:} Update $\rho_{k+1}=\alpha \rho_k$ for a constant $\alpha\in(0,1)$ and repeat until $\varepsilon_{k}=\max(\bar \varepsilon_{1},\bar \varepsilon_{2})>0$, or until no further improvement is observed.
\end{enumerate}
Degrees of \(U,\Lambda,M,S_i\) are fixed a priori; if progress stalls, one can increase degrees or reinitialize \(U_0\). This alternating procedure preserves convexity at each step and typically yields an increase in the slack variable $\varepsilon$. Note that variables $(\alpha_{1:I},\delta,\varepsilon)$ appear linearly in the constraints and therefore can be updated in each step. 

For the update of radius $\rho_k$, we set \(\rho_{k+1}=\alpha\,\rho_k\) with a fixed \(\alpha\in(0,1)\). This geometric schedule monotonically shrinks the disturbance ball \(W_0\) while preserving \(\rho_k>0\) for all \(k\). As \(\rho_k\) decreases, the robust one-step decrease constraint becomes easier to satisfy. Hence, we iterate until the slack \(\varepsilon\) becomes strictly positive. In practice, choosing \(\alpha\) close to \(1\) yields a homotopy with gradual changes that keeps \(W_0\) of nonzero probability (since \(W\) contains an open ball around the origin) while avoiding overly conservative radii.

We transformed the probabilistic variant condition in \textbf{V2} to the tractable optimization in~\eqref{eq:opt:U} by restricting disturbances to a ball of nonzero probability, and encode other  requirements via the S-procedure with SOS multipliers, and resolve the remaining bilinear couplings through an alternating scheme. Each substep reduces to a convex semidefinite optimization problem that can be solved with standard SOS toolchains, yielding a practical, solver-ready pipeline. The algorithm is summarized in~\eqref{alg:variant-sos}.

\begin{algorithm}[t]
\caption{SOS-Based Variant Function}
\label{alg:variant-sos}
\begin{algorithmic}[1]
\State \textbf{Input:} Polynomial dynamics $f(x,w)$, target polynomials $g_i(x)$ for $i\in\{1,\ldots,I\}$, desired degrees for $U,\Lambda,M,S_i$, initial guess $U_0$, initial radius $\rho_0>0$ and a constant $\alpha\in(0,1)$.

\State \textbf{Output:} Variant $U(x)$, multipliers $\Lambda(x,w)$, $M(x,w)$, $S_{1:I}(x)$, radius $\rho^\star>0$ and margins $\delta>0$, $\alpha_{1:I}>0$ satisfying \eqref{eq:U:SOS}–\eqref{eq:GU:SOS}.
\State \textbf{// Step 1: Define polynomial structure}
\State Fix $U\in\mathbb{R}[x]$.
\State \textbf{// Step 2: Initialization}
\State Set $k\gets 0$, $(U,\rho)\gets(U_0,\rho_0)$.
\State \textbf{// Step 3: Alternating optimization}
\Repeat
  \State \textbf{Multiplier step:}
  \State  Solve \eqref{eq:opt:U} for the fixed $(U_k,\rho_k)$ and obtain $$(\Lambda_k,M_k,S_{1:I,k},\delta_{2k-1},\alpha_{1:I,2k-1},\bar\varepsilon_{1}).$$
  \State \textbf{Variant step:}
  \State  Solve \eqref{eq:opt:U} for the fixed $(\Lambda_k,M_k,S_{1:I,k},\rho_k)$ and obtain $$(U_{k+1},\delta_{2k},\alpha_{1:I,2k},\bar \varepsilon_{2}).$$
  \State $\rho_{k+1}\gets \alpha\rho_{k}$ 
  \State $\varepsilon_k\gets \max(\bar \varepsilon_{1},\bar \varepsilon_{2})$
  \State $k\gets k+1$
\Until{$\varepsilon_k>0$ or no further improvement in $\varepsilon$}
\State \textbf{Return:} $(U,\rho^\star)$ if $\varepsilon_k>0$; otherwise increase degrees and reinitialize.
\end{algorithmic}
\end{algorithm}

\begin{Remark}
    Example~\ref{example1} showed that assuming fixed $\delta(r)$ and $\epsilon(r)$ in the variant condition~\textbf{V2} may fail on unbounded state spaces, i.e., the probability of a sufficient one-step decrease can tend to zero as $\|x\|\to\infty$, even when a valid variant exists. However, we used constant $\delta$ and $\epsilon$ in presenting the SOS approach in the paper for clarity and simplicity, and as we will demonstrate in the simulation section, constant margins do succeed for our case studies. On the other hand, the presented SOS approach does not rely only on constant margins, i.e., the algorithms can be modified to also consider $\bar\delta(\cdot)$ and/or $\bar\epsilon(\cdot)$ as being nonnegative polynomials of $x$, design their coefficients, and then get $\delta(r) = \min_x\{\bar\delta(x)\mid V(x)\le r\}$ and $\epsilon(r) = \min_x\{\bar\epsilon(x)\mid V(x)\le r\}$ satisfying the variant condition. This preserves the practicality of the method while aligning with the full generality of \textbf{V2}.
\end{Remark}

\section{Illustrative Example}\label{sec:illustrative}
In this section, we present two illustrative examples that demonstrate the proposed SOS-based framework for certifying almost-sure reachability. Both examples involve two-dimensional stochastic polynomial systems but differ in the way uncertainty influences the system evolution. All SOS programs were modeled in SOSTOOLS~\citep{prajna2004sostools} and solved with SDPT3 in Matlab.

\subsection{Polynomial System with Additive Disturbance}
 Consider the following system with additive disturbance,
\begin{align*}
    x^+_{1} &= 0.3 x_{1} + 0.5 x_{2}^3 + w_{1}, \\
    x^+_{2} &= 0.8 x_{2} + w_{2},
\end{align*}
where \( x = [x_{1}, x_{2}]^\top \in \mathbb{R}^2 \) is the system state and \( w = [w_{1}, w_{2}]^\top \) is a zero-mean i.i.d. disturbance vector uniformly distributed over $[-1,1]^2$. We consider the target set $G=\{x \in \mathbb{R}^2\,|\, x^\top x<1\}$. 

\begin{figure}[t!]
   \centering
    \includegraphics[width=0.4\textwidth]{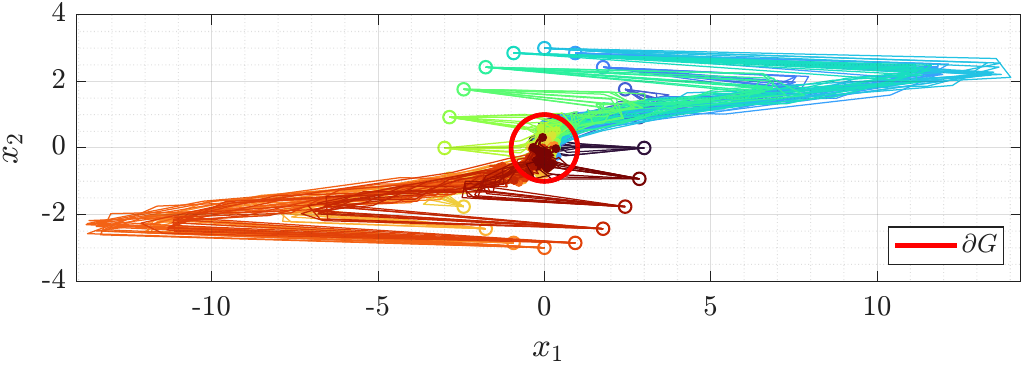}
    \caption{Stochastic trajectories of the system from various initial states.}
    \label{fig:trajectories}
\end{figure}

Figure~\ref{fig:trajectories} shows sample trajectories from multiple initial conditions under i.i.d. noise, illustrating that these paths enter the target set.

\begin{figure}[t!]
    \centering
    \includegraphics[width=0.4\textwidth]{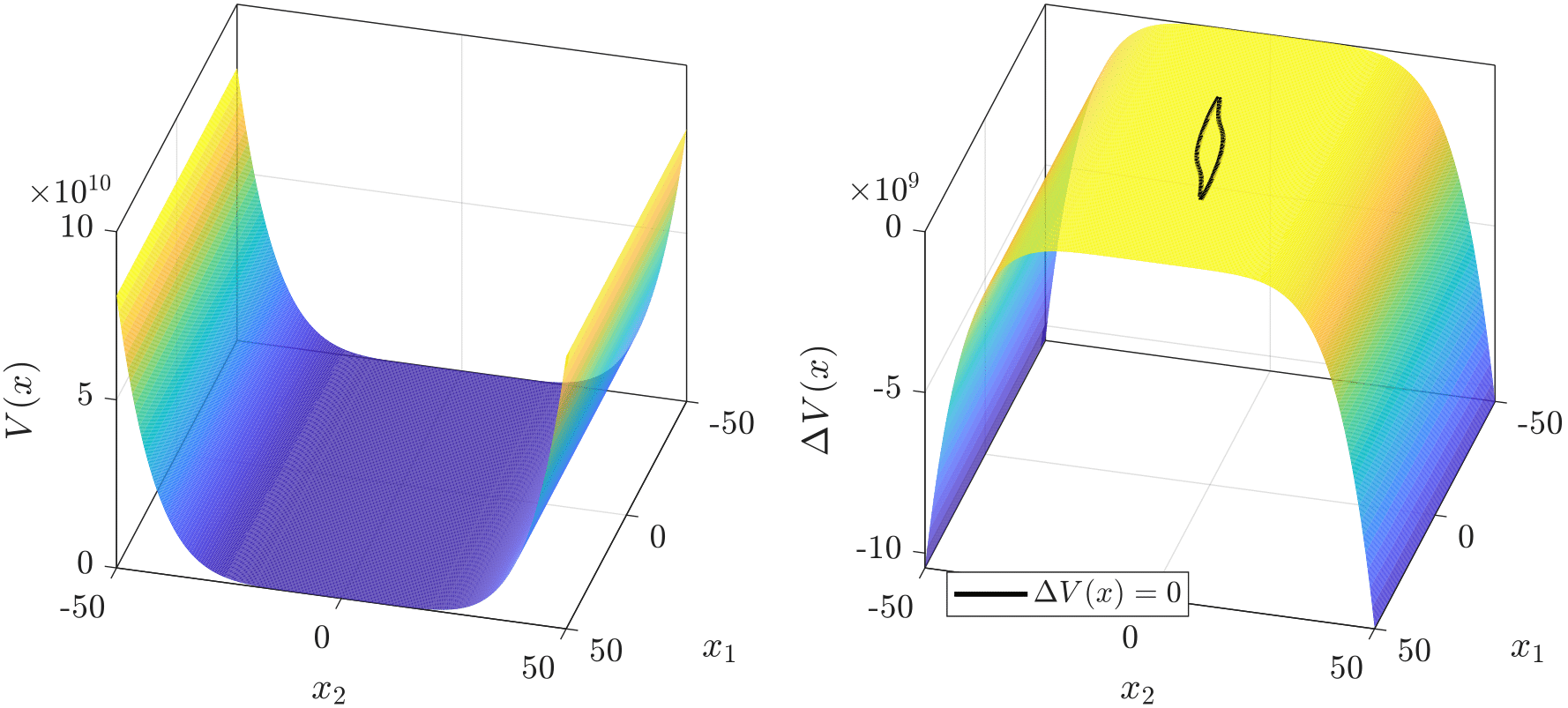}
    \caption{SOS-based drift certification. \textbf{Left:} The certified drift function \( V(x) \). This shows that \( V(x) \) grows radially, satisfying the unboundedness condition. \textbf{Right:} The drift decrement \( \Delta V(x) = \mathbb{E}[V(f(x, w))] - V(x) \). The black contour marks the zero level set and outside of the contour is for $\Delta V(x)\leq 0$.}
    \label{fig:Vx}
\end{figure}

We fix a degree-6 template for \(V\) and solve the drift SOS constraints from Section~\ref{sec:SOS}. The obtained drift is
\begin{align*}
V(x) = & 12.61x_1^2-64.19x_1x_2 +788.27x_2^2
\\ &+20.99x^4_2+5.18x^6_2.
\end{align*}
Note that by considering the quadratic terms, $V(x)$ can be lower bounded by the quadratic function $x^\top A x$ where $A=\begin{bmatrix} 12.61 & -32.095\\ -32.095 & 788.27 \end{bmatrix}$ with $\lambda_{\text{min}}(A)=11.28$. Therefore $V(x)$ is radially unbounded. Figure~\ref{fig:Vx} (left) visualizes the drift function \(V(x)\) and its corresponding decrement \(\Delta V(x)\).  The right figure shows \(\Delta V(x) \leq 0\) and the zero contour of \(\Delta V\) demarcates a compact set outside which the expected change is nonpositive, certifying the drift condition.

We next synthesize a variant \(U\) using Algorithm~\ref{alg:variant-sos}.  For the variant certificate, we use a degree-6 template for \(U\) and degree-$2$ SOS multipliers \(\Lambda(x,w),\,M(x,w),\, S(x)\). The obtained variant $U(x)$ is
\begin{align*}
   &U(x)\!=\! 3.37 x_1 \!-\!\! 1.67 x_2 \!+\! 218.34 x_1^{2}x_2^{2}\! +\! 7.06 x_1^{2}x_2^{3} \!-\! 5.12 x_1^{3}x_2^{2} \\
&\!+\! 41.86 x_1^{2}x_2^{4}\! \!-\! \!20.95 x_1^{3}x_2^{3}\!\! +\! \!3.21 x_1^{4}x_2^{2} \!\!-\!\! 31.02 x_1x_2\! \!- \!\!30.64 x_1x_2^{2} \\
&\!- \!46.03 x_1^{2}x_2 \!\!- \!\!16.5 x_1x_2^{3} \!\!- \!\!167.51 x_1^{3}x_2 \!\!+ \!28.23 x_1x_2^{4}\!\! +\! 0.51 x_1^{4}x_2 \\
&- \!0.04 x_1^{5}x_2 \!+\! 84.24 x_1^{2} \!+ \!12.05 x_1^{3}\! +\! 4.41 x_2^{2}\! +\! 27.83 x_1^{4} \\
&+ \!39.48 x_2^{3}\! -\! 0.01 x_1^{5}\! +\! 651.48 x_2^{4}\! - \!0.02 x_1^{6} \!-\! 0.13 x_2^{5}\!-\! 11.64.
\end{align*}

\begin{figure}[t!]
    \centering
    \includegraphics[width=0.4\textwidth]{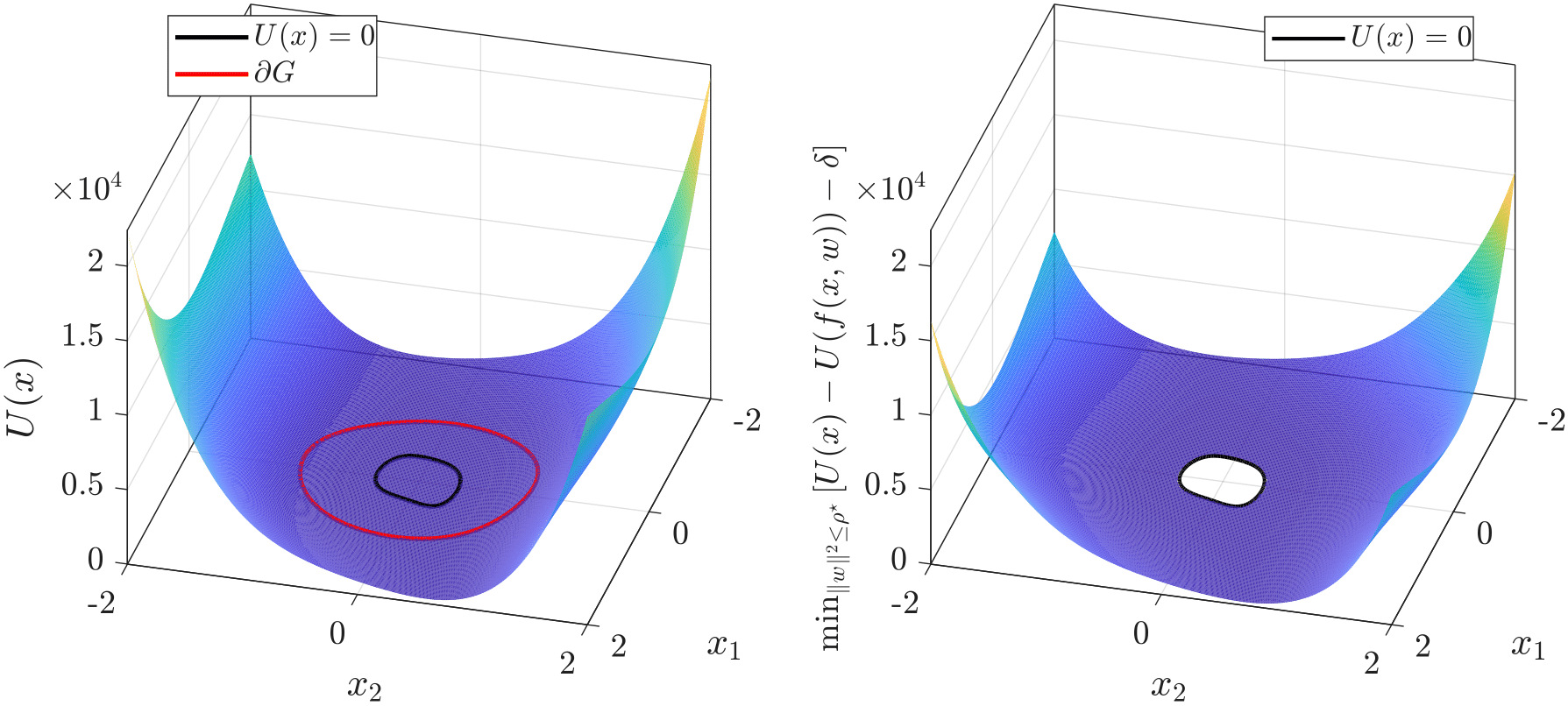}
\caption{SOS-based variant certification. \textbf{Left:} Surface of $U(x)$ with the contour $U(x)=0$ (black) and the target-set boundary $\partial G$ (red), showing $\{x:\,U(x)\le 0\}\subseteq G$. \textbf{Right:} $\min_{\|w\|^2\le \rho^\star}\!\big[\,U(x)-U(f(x,w))-\delta\,\big]$ over the state space. For all states with $U(x)>0$, the quantity $U(x)-U(f(x,w))-\delta$ is nonnegative for every disturbance $\|w\|^2\le \rho^\star$, certifying a uniform one-step decrease on a non-empty set.}
\label{fig:variant}
\end{figure}

\begin{figure}[t!]
    \centering
    \includegraphics[width=0.4\textwidth]{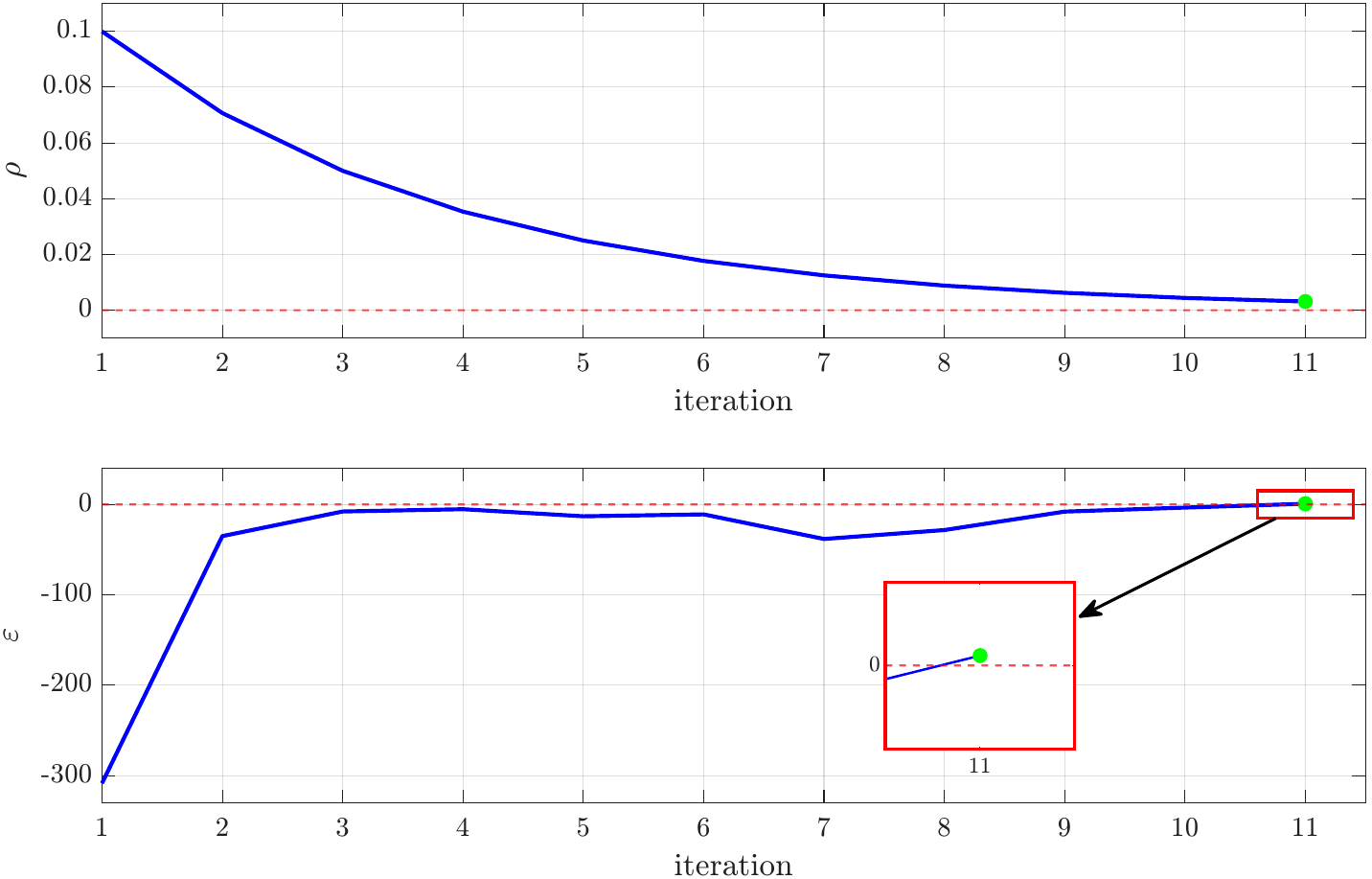}
   \caption{Evolution of $\rho_k$ (top) and the slack $\varepsilon_k$ (bottom) across iterations of Algorithm~\ref{alg:variant-sos}.}
    \label{fig4}
\end{figure}

Figure~\ref{fig:variant} (left) shows the computed variant $U(x)$ together with the target-set boundary $\partial G$ and the contour $U(x)=0$. The plot indicates that the target set contains all states with $U(x)\le 0$, i.e., $\{x:\, U(x)\le 0\}\subset G$. Figure~\ref{fig:variant} (right) displays the map $x\mapsto \min_{\|w\|^2\le \rho^\star}\big[\,U(x)-U(f(x,w))-\delta\,\big]$ over the state space. It can be seen that for all $x$ with $U(x)>0$, we have $U(x)-U(f(x,w))-\delta\ge 0$ for all $w$ satisfying $\|w\|^2\le \rho^\star$. Figure~\ref{fig4} tracks the evolution of the disturbance radius $\rho_k$ (top) and the slack $\varepsilon_k$ (bottom) across iterations of Algorithm~\ref{alg:variant-sos}. The bottom figure shows that $\varepsilon_k$ becomes strictly positive at the final iteration, certifying feasibility, while the top figure shows a decrease in $\rho_k$ that remains positive at termination (i.e., both $\rho_k>0$ and $\varepsilon_k>0$ at the last step).

\subsection{Polynomial System with Multiplicative  Disturbance}
Consider now the following stochastic nonlinear system with multiplicative noise in both state coordinates,
\begin{align*}
    x^+_{1} &= x_{1} + 0.05\big(-1.1x_{1} + 0.15x_{1}x_{2} - 0.005x_{1}^3 + w_{1}x_{1}\big), \\
    x^+_{2} &= x_{2} + 0.05\big(-0.9x_{2} + 0.12x_{1}^2 - 0.006x_{2}^3 + w_{2}x_{2}\big),
\end{align*}
where \(x = [x_{1},x_{2}]^\top \in \mathbb{R}^2\) is the state vector and \(w = [w_{1},w_{2}]^\top\) is a zero-mean i.i.d.\ disturbance vector uniformly distributed over $[-0.5,0.5]^2$. We consider the target set \(G = \{x \in \mathbb{R}^2 \mid x^\top x < 1\}\). 

\begin{figure}[t!]
    \centering
    \includegraphics[width=0.4\textwidth]{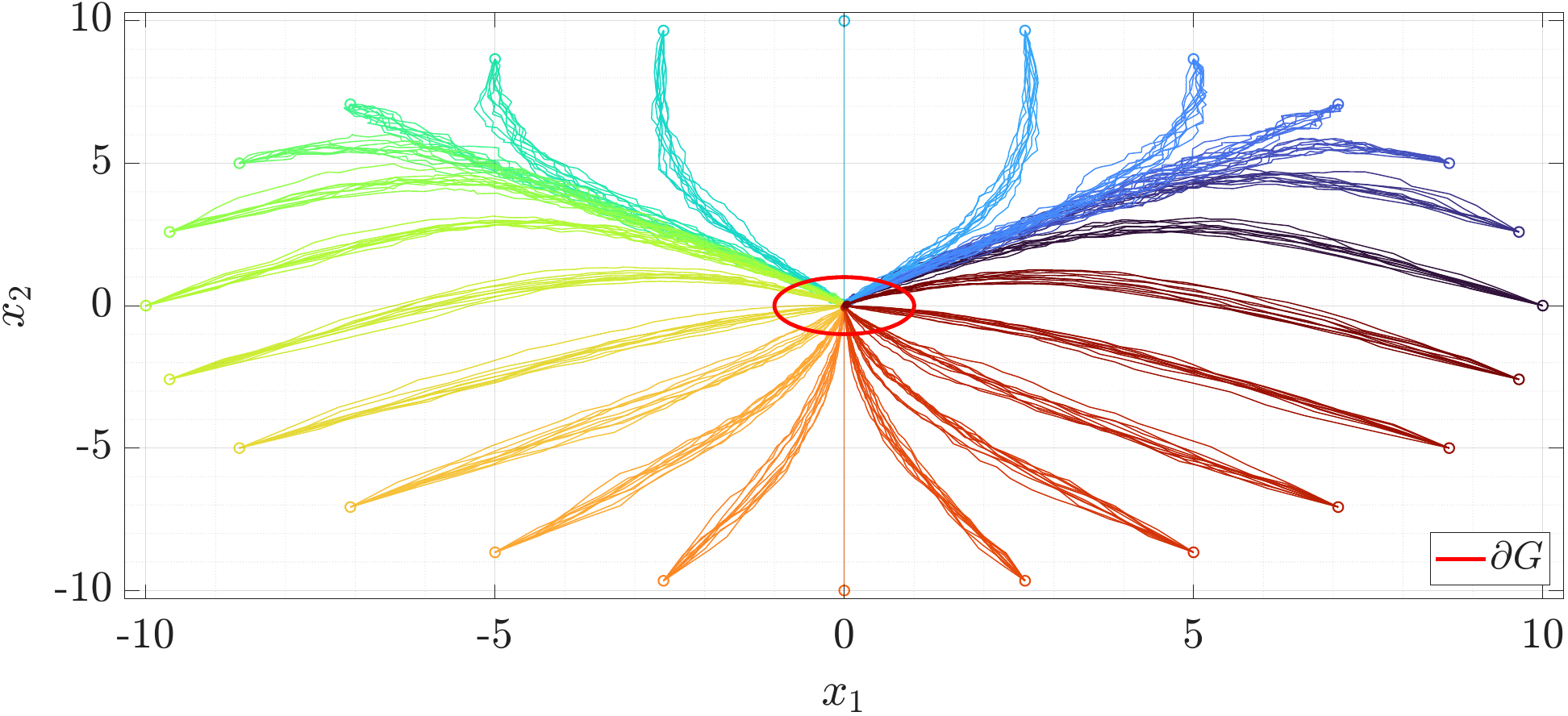}
    \caption{Stochastic trajectories of the multiplicative-noise system from various initial conditions. Each color corresponds to a distinct starting point on a circle of radius~$10$.}
    \label{fig:mult_traj}
\end{figure}

Figure~\ref{fig:mult_traj} illustrates stochastic trajectories generated from initial points uniformly distributed on a circle of radius~$10$. Each colored curve represents an independent realization, and it can be observed that all sample paths eventually enter the target set despite the state-dependent stochastic fluctuations.

\begin{figure}[t!]
    \centering
    \includegraphics[width=0.4\textwidth]{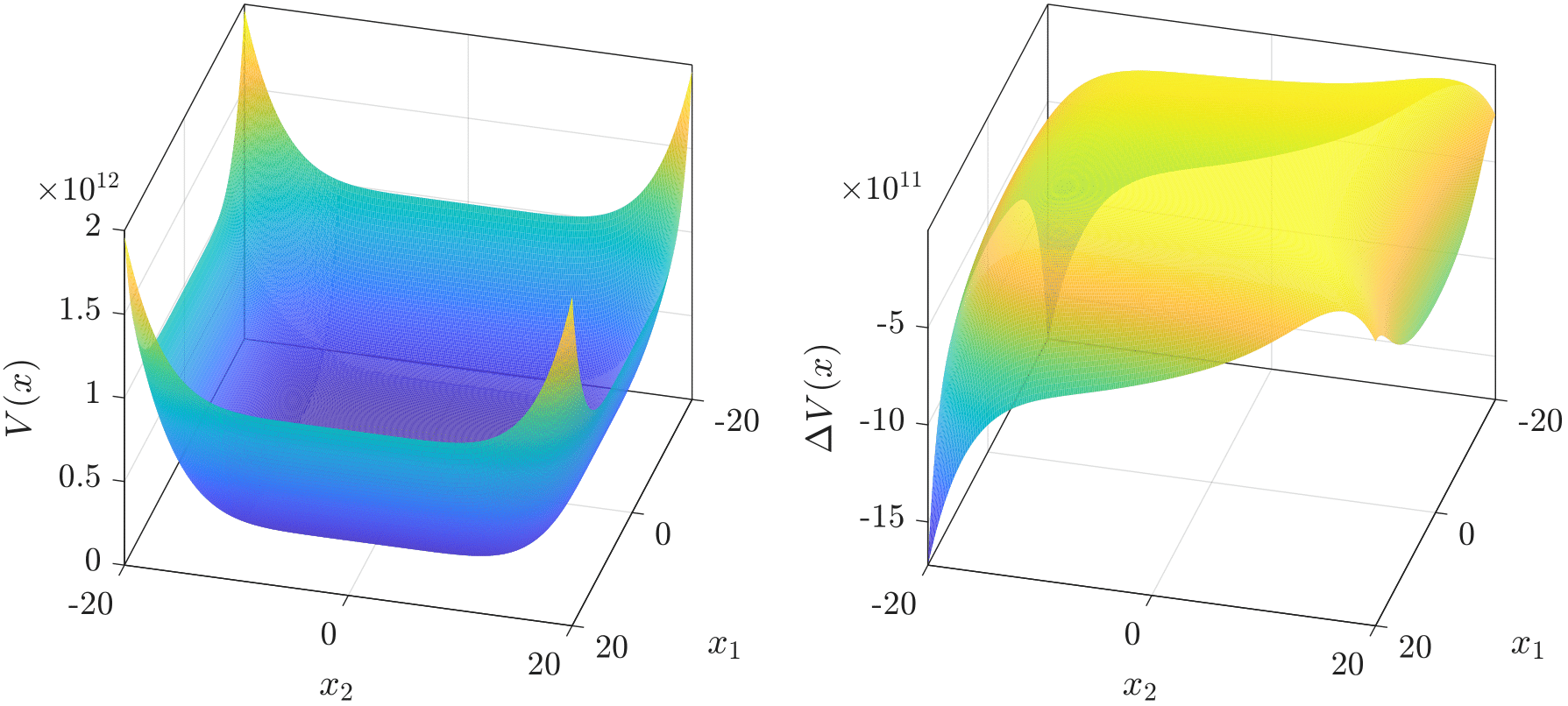}
    \caption{SOS-based drift certification for the multiplicative-noise system. \textbf{Left:} The computed drift function \(V(x)\) grows radially, ensuring unboundedness. \textbf{Right:} The drift decrement \(\Delta V(x)\), verifying the expected non-increase of \(V(x)\) outside it.}
    \label{fig:mult_V}
\end{figure}

For the drift function, we fix a 6-degree  polynomial template and solve the SOS constraints enforcing radial unboundedness and non-increasing expectation outside a compact set. The resulting drift function is of the form
\begin{align*}
V(x)\! =\!{}& 15425.40x_1^6 \!+ \!1221.10x_1^4 \!+\! 933.08x_1^2x_2^2 \!+\! 856.16x_1^2 \\
& -\! 0.05x_1x_2\! +\! 15281.28x_2^6 \!+\! 1201.01x_2^4 \!+\! 843.80x_2^2,
\end{align*}
and is radially unbounded, satisfying the drift condition in~\eqref{eq:drift}.  Figure~\ref{fig:mult_V} visualizes the obtained drift certificate. The left plot shows the radially growing surface of $V(x)$, while the right plot depicts the drift decrement $\Delta V(x) = \mathbb{E}[V(f(x,w))]-V(x)$.

We then synthesize a {variant} function \(U\) using the alternating SOS optimization. The resulting polynomial variant function is
\begin{align*}
&U(x) \!=\!1.56x_1^6 \!+\! 1.74x_1^4x_2^2 \!- \!3.59x_1^4x_2 \!+\! 5.22x_1^4\! +\! 1.82x_1^2x_2^4 \\
& -\! 2.85x_1^2x_2^3\! +\! 2.65x_1^2x_2^2\! -\! 0.99x_1^2x_2\! +\! 19.85x_1^2\!+ \!1.20x_2^6  \\
& - \!0.24x_2^5 \!+\! 1.85x_2^4 \!-\! 0.31x_2^3\! +\! 28.84x_2^2 \!+\! 0.02x_2 \!-\! 2.94,
\end{align*}
and satisfies the inclusion $\{x\mid U(x)\le0\}\subset G$ as well as the one-step robust decrease property. 

\begin{figure}[t!]
    \centering
    \includegraphics[width=0.4\textwidth]{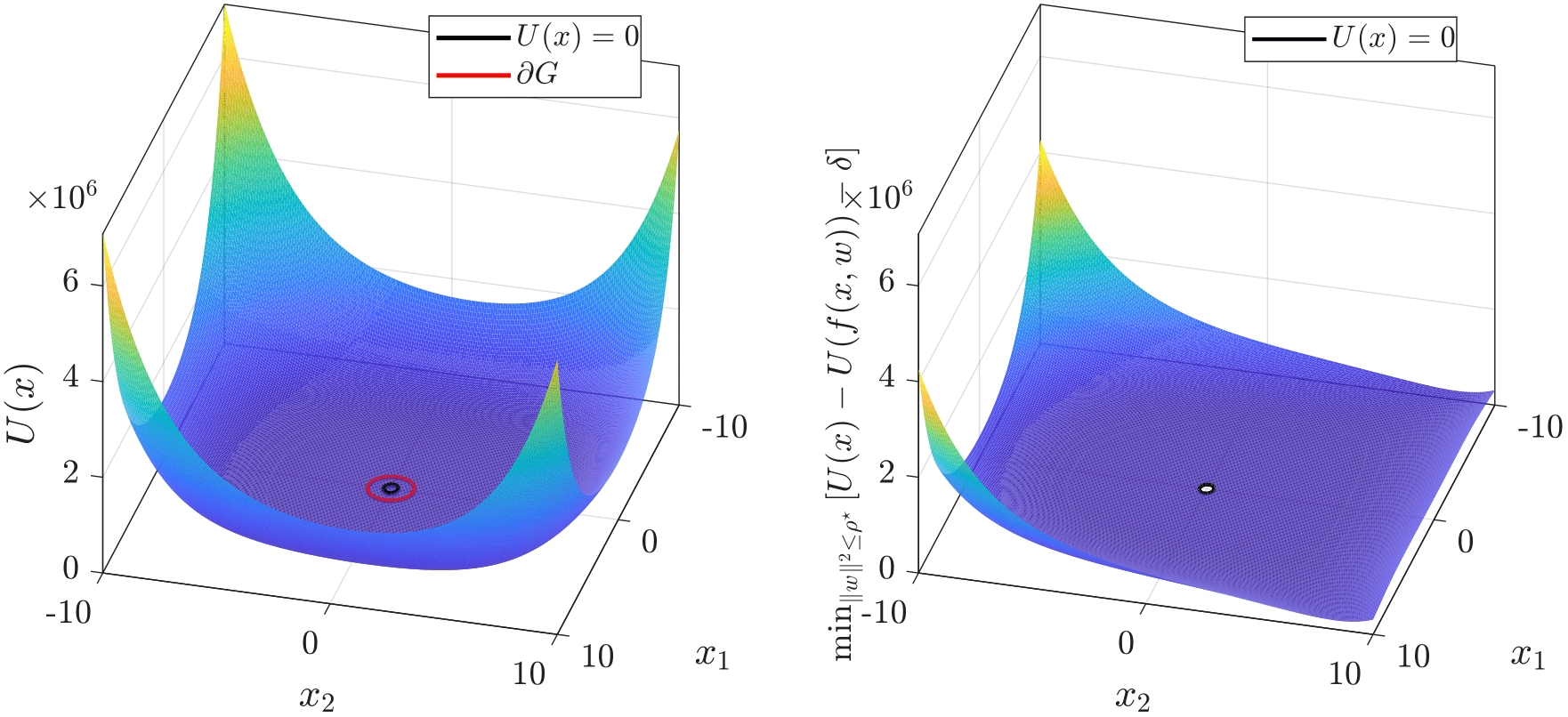}
    \caption{SOS-based variant certification for the multiplicative-noise system. \textbf{Left:} Surface of $U(x)$ with the contour $U(x)=0$ (black) and target boundary $\partial G$ (red), confirming $\{x:U(x)\le0\}\subset G$. \textbf{Right:} The minimum robust one-step decrease $\min_{\|w\|^2\le\rho^\star}[U(x)-U(f(x,w))-\delta]$ showing nonnegativity over $U(x)>0$, certifying almost-sure reachability under multiplicative noise.}
    \label{fig:mult_U}
\end{figure}

Figure~\ref{fig:mult_U} (left) presents the computed surface $U(x)$, where the contour $U(x)=0$ (black) lies strictly inside the target boundary $\partial G$ (red). This confirms that $\{U(x)\le0\}\subset G$. The right plot shows that for all $x$ with $U(x)>0$, $U(x)-U(f(x,w))-\delta$ remains nonnegative across all admissible disturbances, certifying a uniform one-step decrease on a subset of disturbances with positive probability.

\begin{figure}[t!]
    \centering
    \includegraphics[width=0.4\textwidth]{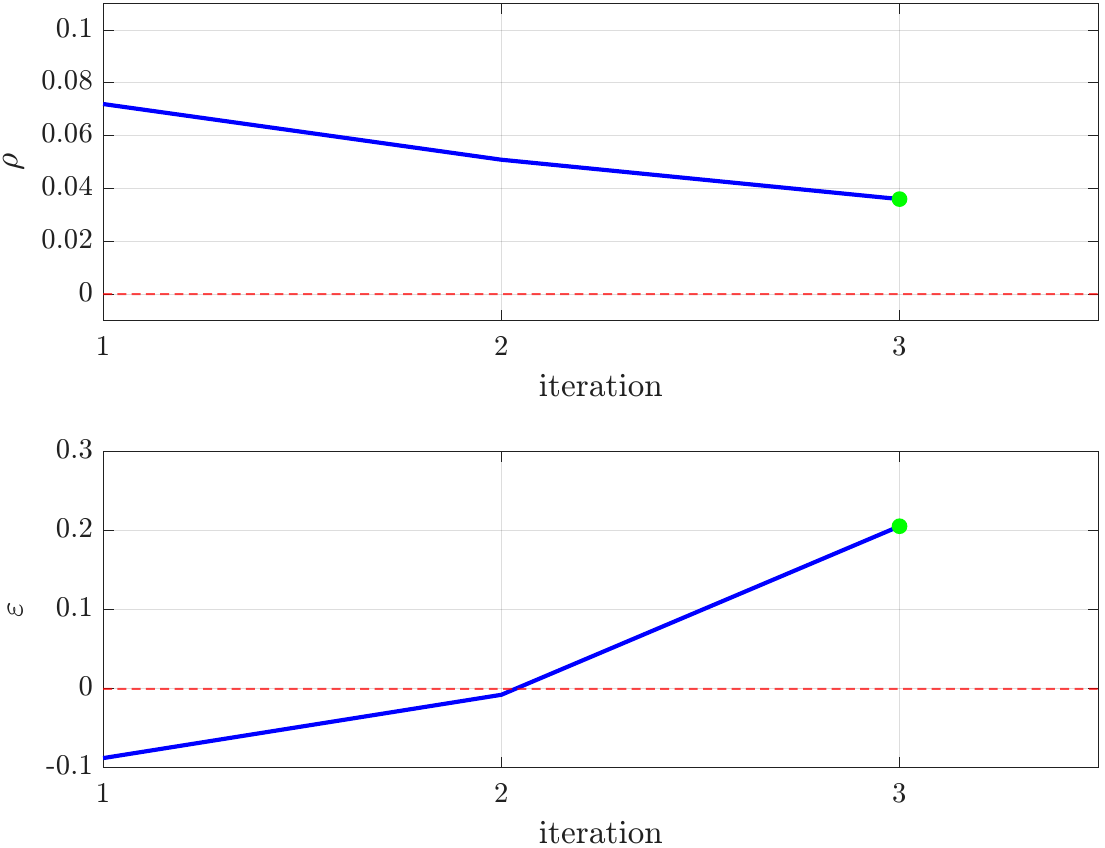}
    \caption{Evolution of the disturbance radius $\rho_k$ (top) and slack variable $\varepsilon_k$ (bottom) during the alternating SOS synthesis. The final positive $\varepsilon_k$ certifies feasibility of the variant condition for the multiplicative-noise dynamics.}
    \label{fig:mult_rho_eps}
\end{figure}

Figure~\ref{fig:mult_rho_eps} depicts the evolution of the disturbance radius $\rho_k$ and the slack variable $\varepsilon_k$ over successive iterations of the alternating SOS procedure. The radius $\rho_k$ decreases geometrically while remaining positive, and $\varepsilon_k$ becomes strictly positive at the final iteration, confirming successful certificate synthesis.


\section{Discussion and Conclusion}
\label{sec:conclusion}
We presented a Sum-of-Squares (SOS) approach for discrete-time polynomial stochastic systems that certifies almost-sure reachability by casting the drift–variant conditions as tractable semidefinite programs. Once polynomial templates and degrees are fixed for the drift and variant functions, the constraints are encoded via an S-procedure, yielding a solver-ready formulation. Bilinear couplings among decision variables are handled with an alternating scheme. 

While the provided computational approach makes almost-sure reachability certificates practical, it has its own limitations. For example, reliance on polynomial templates and the S-procedure may be conservative, scalability degrades with polynomial degree and state dimension, and the SOS optimizations may be numerically unstable \citep{roux2018validating}. Therefore, richer certificate parameterizations and nonpolynomial templates (e.g., neural networks, rational, exponential, or logarithmic functions) can be selected for the candidate certificates. Note that for such parametrizations, the techniques described in this paper can be used for handling the existential quantifier on the compact set in the drift condition and for eliminating the probability operator in the variant condition. Therefore, our results provide a basis also for applying computational methods that are based on counter-example guided inductive synthesis and neural templates \citep{abate2021fossil,nejati2020compositional,abate2024safe}. Our future investigation also includes the control synthesis problem for ensuring almost sure reachability, the existence of such control policies under appropriate assumptions, and deriving such certificates for continuous-time systems without the need for employing time discretization techniques that do not generally preserve temporal properties of the system.

\bibliography{SOS}

\begin{thebibliography}{24}
\providecommand{\natexlab}[1]{#1}
\providecommand{\url}[1]{\texttt{#1}}
\providecommand{\urlprefix}{URL }
\expandafter\ifx\csname urlstyle\endcsname\relax
  \providecommand{\doi}[1]{doi:\discretionary{}{}{}#1}\else
  \providecommand{\doi}{doi:\discretionary{}{}{}\begingroup \urlstyle{rm}\Url}\fi

\bibitem[{Abate et~al.(2021)Abate, Ahmed, Edwards, Giacobbe, and Peruffo}]{abate2021fossil}
Abate, A., Ahmed, D., Edwards, A., Giacobbe, M., and Peruffo, A. (2021).
\newblock Fossil: a software tool for the formal synthesis of {L}yapunov functions and barrier certificates using neural networks.
\newblock In \emph{Proceedings of the 24th international conference on hybrid systems: computation and control}, 1--11.

\bibitem[{Abate et~al.(2024)Abate, Bogomolov, Edwards, Potomkin, Soudjani, and Zuliani}]{abate2024safe}
Abate, A., Bogomolov, S., Edwards, A., Potomkin, K., Soudjani, S., and Zuliani, P. (2024).
\newblock Safe reach set computation via neural barrier certificates.
\newblock \emph{IFAC-PapersOnLine}, 58(11), 107--114.

\bibitem[{Ahmadi et~al.(2011)Ahmadi, Krstic, and Parrilo}]{ahmadi2011globally}
Ahmadi, A.A., Krstic, M., and Parrilo, P.A. (2011).
\newblock A globally asymptotically stable polynomial vector field with no polynomial {L}yapunov function.
\newblock In \emph{2011 50th IEEE Conference on Decision and Control and European Control Conference}, 7579--7580. IEEE.

\bibitem[{Ahmadi and Parrilo(2012)}]{ahmadi2012convex}
Ahmadi, A.A. and Parrilo, P.A. (2012).
\newblock A convex polynomial that is not sos-convex.
\newblock \emph{Mathematical Programming}, 135(1), 275--292.

\bibitem[{Ames et~al.(2019)Ames, Coogan, Egerstedt, Notomista, Sreenath, and Tabuada}]{ames2019control}
Ames, A.D., Coogan, S., Egerstedt, M., Notomista, G., Sreenath, K., and Tabuada, P. (2019).
\newblock Control barrier functions: Theory and applications.
\newblock In \emph{2019 18th European control conference (ECC)}, 3420--3431. Ieee.

\bibitem[{Baier and Katoen(2008)}]{baier2008principles}
Baier, C. and Katoen, J.P. (2008).
\newblock \emph{Principles of model checking}.
\newblock MIT press.

\bibitem[{Jagtap et~al.(2020)Jagtap, Soudjani, and Zamani}]{jagtap2020formal}
Jagtap, P., Soudjani, S., and Zamani, M. (2020).
\newblock Formal synthesis of stochastic systems via control barrier certificates.
\newblock \emph{IEEE Transactions on Automatic Control}, 66(7), 3097--3110.

\bibitem[{Junges et~al.(2021)Junges, Jansen, and Seshia}]{junges2021enforcing}
Junges, S., Jansen, N., and Seshia, S.A. (2021).
\newblock Enforcing almost-sure reachability in {POMDPs}.
\newblock In \emph{International Conference on Computer Aided Verification}, 602--625. Springer.

\bibitem[{Kordabad et~al.(2024)Kordabad, Charitidou, Dimarogonas, and Soudjani}]{kordabad2024control}
Kordabad, A.B., Charitidou, M., Dimarogonas, D.V., and Soudjani, S. (2024).
\newblock Control barrier functions for stochastic systems under signal temporal logic tasks.
\newblock In \emph{2024 European Control Conference (ECC)}, 3213--3219. IEEE.

\bibitem[{Kordabad et~al.(2025)Kordabad, Majumdar, Motwani, and Soudjani}]{kordabad2025certificates}
Kordabad, A.B., Majumdar, R., Motwani, H.J., and Soudjani, S. (2025).
\newblock On certificates for almost sure reachability in stochastic systems.
\newblock \emph{arXiv preprint arXiv:2507.20194}.

\bibitem[{Lavaei et~al.(2022)Lavaei, Soudjani, Abate, and Zamani}]{lavaei2022automated}
Lavaei, A., Soudjani, S., Abate, A., and Zamani, M. (2022).
\newblock Automated verification and synthesis of stochastic hybrid systems: A survey.
\newblock \emph{Automatica}, 146, 110617.

\bibitem[{Lin et~al.(2023)Lin, Chong, and Murguia}]{lin2023secondary}
Lin, Y., Chong, M.S., and Murguia, C. (2023).
\newblock Secondary controller design for the safety of nonlinear systems via sum-of-squares programming.
\newblock In \emph{2023 62nd IEEE Conference on Decision and Control (CDC)}, 7649--7654. IEEE.

\bibitem[{Lin et~al.(2025)Lin, Chong, and Murguia}]{lin2025modified}
Lin, Y., Chong, M.S., and Murguia, C. (2025).
\newblock Modified control barrier function for quadratic program based control design via sum-of-squares programming.
\newblock \emph{arXiv preprint arXiv:2504.19796}.

\bibitem[{Majumdar and Sathiyanarayana(2025)}]{majumdar2025sound}
Majumdar, R. and Sathiyanarayana, V. (2025).
\newblock Sound and complete proof rules for probabilistic termination.
\newblock \emph{Proceedings of the ACM on Programming Languages}, 9(POPL), 1871--1902.

\bibitem[{Majumdar et~al.(2024)Majumdar, Sathiyanarayana, and Soudjani}]{majumdar2024necessary}
Majumdar, R., Sathiyanarayana, V., and Soudjani, S. (2024).
\newblock Necessary and sufficient certificates for almost sure reachability.
\newblock \emph{IEEE Control Systems Letters}.

\bibitem[{Meyn and Tweedie(2012)}]{meyn2012markov}
Meyn, S.P. and Tweedie, R.L. (2012).
\newblock \emph{Markov chains and stochastic stability}.
\newblock Springer Science \& Business Media.

\bibitem[{Nejati et~al.(2020)Nejati, Soudjani, and Zamani}]{nejati2020compositional}
Nejati, A., Soudjani, S., and Zamani, M. (2020).
\newblock Compositional construction of control barrier certificates for large-scale stochastic switched systems.
\newblock \emph{IEEE Control Systems Letters}, 4(4), 845--850.

\bibitem[{Papachristodoulou and Prajna(2005)}]{papachristodoulou2005tutorial}
Papachristodoulou, A. and Prajna, S. (2005).
\newblock A tutorial on sum of squares techniques for systems analysis.
\newblock In \emph{Proceedings of the 2005, American Control Conference, 2005.}, 2686--2700. IEEE.

\bibitem[{Parrilo(2000)}]{parrilo2000structured}
Parrilo, P.A. (2000).
\newblock \emph{Structured semidefinite programs and semialgebraic geometry methods in robustness and optimization}.
\newblock California Institute of Technology.

\bibitem[{Prajna(2004)}]{prajna2004sostools}
Prajna, S. (2004).
\newblock {SOSTOOLS}: {S}um of squares optimization toolbox for {MATLAB}.
\newblock \emph{http://www. mit. edu/\~{} parrilo/sostools/index. html}.

\bibitem[{Roux et~al.(2018)Roux, Voronin, and Sankaranarayanan}]{roux2018validating}
Roux, P., Voronin, Y.L., and Sankaranarayanan, S. (2018).
\newblock Validating numerical semidefinite programming solvers for polynomial invariants.
\newblock \emph{Formal Methods in System Design}, 53(2), 286--312.

\bibitem[{Schneeberger et~al.(2023)Schneeberger, D{\"o}rfler, and Mastellone}]{schneeberger2023sos}
Schneeberger, M., D{\"o}rfler, F., and Mastellone, S. (2023).
\newblock {SOS} construction of compatible control {L}yapunov and barrier functions.
\newblock \emph{IFAC-PapersOnLine}, 56(2), 10428--10434.

\bibitem[{Schweidel et~al.(2022)Schweidel, Yin, Smith, and Arcak}]{safe2022}
Schweidel, K.S., Yin, H., Smith, S.W., and Arcak, M. (2022).
\newblock Safe-by-design planner–tracker synthesis with a hierarchy of system models.
\newblock \emph{Annual Reviews in Control}, 53, 138--146.

\bibitem[{Yin et~al.(2021)Yin, Arcak, Packard, and Seiler}]{yin2021backward}
Yin, H., Arcak, M., Packard, A., and Seiler, P. (2021).
\newblock Backward reachability for polynomial systems on a finite horizon.
\newblock \emph{IEEE Transactions on Automatic Control}, 66(12), 6025--6032.

\end{thebibliography}
\end{document}